\documentclass[oneside,notitlepage,12pt]{article}

\usepackage{amssymb}
\usepackage[leqno]{amsmath}
\usepackage{amsfonts}
\usepackage{amsopn}
\usepackage{amstext}
\usepackage{amsthm}
\usepackage{enumitem}

\usepackage{stmaryrd}

\usepackage{tikz}
\usetikzlibrary{cd}
\providecommand{\ar}{\arrow}

\usepackage[colorlinks, backref]{hyperref}

\usepackage{calrsfs}
\usepackage{fourier-orns}
\usepackage{hieroglf} 
\usepackage{clock} 

\renewcommand{\Bbb}{\mathbb}

\newenvironment{pf}{\begin{proof}}{\end{proof}}

\newcommand{\Zee}{{\Bbb{Z}}}

\newcommand{\Nat}{{\Bbb{N}}}

\newcommand{\al}{\alpha}

\renewcommand{\phi}{\varphi}
\renewcommand{\rho}{\varrho}

\newcommand{\rest}{\restriction}

\newcommand{\ntr}{{n\in\omega}}

\newcommand{\loe}{\leqslant}

\newcommand{\subs}{\subseteq}

\newcommand{\nnempty}{\ne\emptyset}
\newcommand{\argum}{\:\cdot\:}

\newcommand{\id}[1]{{\operatorname{i\!d}_{#1}}} 

\newcommand{\oraz}{\qquad\text{and}\qquad}

\newcommand{\setof}[2]{\{#1\colon #2\}}

\newcommand{\sett}[2]{\{#1\}_{#2}}
\newcommand{\sn}[1]{\{#1\}} 
\newcommand{\pair}[2]{\langle #1, #2 \rangle} 
\newcommand{\map}[3]{#1\colon #2 \to #3} 
\newcommand{\img}[2]{#1[#2]} 

\newcommand{\M}{\mathbb M}

\providecommand{\nat}{\omega}

\newcommand{\ciag}[1]{{\sett{{#1}_n}{\ntr}}}

\newcommand{\G}{{\mathbb G}}

\newcommand{\End}{\operatorname{End}}

\newcommand{\fK}{{\mathfrak{K}}}

\newcommand{\fC}{{\mathfrak{C}}}

\newcommand{\cmp}{\circ} 

\newcommand{\ob}[1]{\operatorname{Obj}\left(#1\right)}

\newcommand{\acts}{\curvearrowright}
\newcommand{\action}[3]{{#1}^{\acts}_{#2}{#3}} 

\newcommand{\proto}[1]{{\mathbb S_\kappa}}

\newtheorem{tw}{Theorem}[section]
\newtheorem{wn}[tw]{Corollary}
\newtheorem{lm}[tw]{Lemma}
\newtheorem{prop}[tw]{Proposition}

\theoremstyle{definition}
\newtheorem{df}[tw]{Definition}
\newtheorem{ex}[tw]{Example}

\theoremstyle{remark}
\newtheorem{uwgi}[tw]{Remark}

\newcommand{\uv}[1]{\widehat{#1}} 
\newcommand{\un}{\upsilon} 
\newcommand{\uno}{\sigma}

\newcommand{\frg}[1]{|#1|} 
\newcommand{\frgc}{\frg{\cdot}} 
\newcommand{\F}{\mathbf{F}} 

\newcommand{\Ens}{\mathfrak E\! \mathfrak n\! \mathfrak s} 

\newcommand{\nepis}[1]{{#1}^{\heartsuit}}

\newcommand{\Gr}[1]{\G{#1}} 
\newcommand{\MF}{\action{\M}{}{\F}}

\title{Universal monoid actions: the power of freedom}
\author{
{\sc Wies{\l}aw Kubi\'s}\footnote{Research supported by GA \v{C}R grant 20-22230L (Czech Science Foundation).}\\ 
{\small Institute of Mathematics, Czech Academy of Sciences}
}
\date{\today\ \clocktime}

\begin{document}

\maketitle

\begin{abstract}
	Motivated by a recent work of Balcerzak and Kania [Proc. Amer. Math. Soc. 151 (2023) 3737--3742], we show that every countable monoid has a universal action on the free object over a countable infinite set. This is a general result concerning concrete categories with a left adjoint (free) functor. On the way, we introduce an abstract concept of ``being generated by a set''. At the same time we obtain a simpler proof of the result of Darji and Matheron concerning a surjectively universal operator on the classical Banach space $\ell_1$ of summable sequences.
	
	\ 



\end{abstract}

\tableofcontents

\section{Introduction}

The aim of this note is to present a story of universal monoid actions on concrete mathematical objects, in the sense that every action lifts to a single one, naturally defined. This universal action can be viewed as a free construction, as it indeed comes from a left adjoint to the functor that forgets the monoid action.


Let us start with the following recent result, due to Darji and Matheron~\cite{DaMa}:

\begin{quote}
	There exists a non-expansive linear operator $\map U {\ell_1}{\ell_1}$ such that for every separable Banach space $E$, for every non-expansive operator $\map f E E$ there is a non-expansive surjective linear operator $\map q {\ell_1}E$ satisfying $f \cmp q = q \cmp U$.
\end{quote}

For simplicity, we shall say that $f$ \emph{lifts} to $U$ if there is $q$ as above, namely, satisfying $f \cmp q = q \cmp U$.
Given a nonempty set $S$, recall that $\ell_1(S)$ is the Banach space of all summable mappings from $S$ to the real or complex numbers (depending on our choice). A mapping $\map v S \mathbb C$ is \emph{summable}, if $\sum_{s \in S}|v(s)| < +\infty$.
Given $s\in S$, we denote by $\uv s$ the basic vector of $\ell_1(S)$ induced by $s$, namely, the characteristic function of $\sn s$.

We now show how to obtain a surjectively universal operator, using the fact that $\ell_1$ is the free Banach space over the countable infinite set.
The purpose of this note is to extend this result to monoid actions on arbitrary concrete categories.

Let us call an operator $\map{T}{\ell_1(S)}{\ell_1(S)}$ \emph{basic} if it is induced by a self-map $\map \phi S S$ in the sense that $T \uv s = \uv{\phi(s)}$ for every $s \in S$. In this case we shall write $T = \ell_1(\phi)$. Indeed, $\ell_1(\argum)$ can be viewed as a functor from the category of sets into the category of Banach spaces.
Clearly, a basic operator has necessarily norm one and it is a linear isometric embedding whenever $\phi$ is an injection.

\begin{lm}
	Let $E$ be a separable Banach space. Then every non-expansive linear operator $\map f E E$ lifts to a basic operator on $\ell_1$.
\end{lm}

\begin{pf}
	Fix a countable $f$-invariant set $S$ that is dense in the unit ball of $E$. There is a unique linear operator $\map q{\ell_1(S)}E$ such that $q(\uv s) = s$ for every $s \in S$.
	Note that $q$ is surjective, because the image of the unit ball of $\ell_1(S)$ is dense in the unit ball of $E$---this is a consequence of a well-known fact in Banach space theory, see e.g.~\cite[Lemma 2.24]{Fab}.
	Let $\map{T}{\ell_1(S)}{\ell_1(S)}$ be such that $T \uv s = \uv{f(s)}$ for $s\in S$. Then
	$$f(q(\uv s)) = f(s) = q(T \uv s)$$
	for every $s\in S$, therefore $f \cmp q = q \cmp T$. Finally, we may replace $S$ by $\nat$, so that $\ell_1(S) = \ell_1$.
\end{pf}

We are now in a position where the existence of a surjectively universal operator $U$ can be reduced to proving the existence of a surjectively universal self-mapping of a countable set. Namely, knowing that $\ell_1(\argum)$ is a functor, we have the following simple fact.

\begin{lm}
	Let $S$, $T$ be nonempty sets, $\map f S S$, $\map g T T$ and $\map p T S$ be such that $p$ is a surjection and $p \cmp g = f \cmp p$. Then the diagram
	$$\begin{tikzcd}
		\ell_1(T) \ar[rr, "\ell_1(g)"] \ar[d, two heads, "\ell_1(p)"'] & & \ell_1(T) \ar[d, two heads, "\ell_1(p)"] \\
		\ell_1(S) \ar[rr, "\ell_1(f)"]& & \ell_1(S)
,	\end{tikzcd}$$
	is commutative and $\ell_1(p)$ is a norm one projection.
\end{lm}

In other words, under the assumptions above, $\ell_1(f)$ lifts to $\ell_1(g)$,


Let $\map{\un}{\nat^2}{\nat^2}$ be defined by
$\un \pair m n = \pair{m+1}{n}$.
Clearly, $\un$ is one-to-one and all of its orbits are infinite. We will see later that $\un$ is a special case of universal monoid actions.

\begin{lm}\label{LMefivhid}
	$\un$ is surjectively universal. Namely, given a nonempty countable set $S$, given a mapping $\map f SS$, there is a surjection $\map q {\nat^2}S$ such that $q \cmp \un = f \cmp q$.
\end{lm}

\begin{pf}
	Let $S = \ciag s$ and define
	$$q \pair m n = f^m(s_n), \qquad \pair m n \in \nat^2.$$
	We have $q(\un \pair m n) = q \pair {m+1}n = f^{m+1}(s_n)$
	and $f(q \pair m n) = f(f^m(s_n)) = f^{m+1}(s_n)$, therefore $q \cmp \un = f \cmp q$.	
\end{pf}

\begin{wn}\label{WNunivElone}
	The operator $\ell_1(\un)$ is surjectively universal in the category of separable Banach spaces.
\end{wn}

\begin{uwgi}\label{UWGAfgsdhgs}
	An analogue of Lemma~\ref{LMefivhid} and the above corollary holds for uncountable densities, where we replace $\un$ by $\map{\un_\kappa}{\nat \times \kappa}{\nat \times \kappa}$ defined by
	$\un_\kappa \pair n \al = \pair {n+1} \al$. Given a nonempty set $S$ of cardinality $\loe \kappa$, we enumerate it as $\sett{s_\al}{\al<\kappa}$ and the proof of Lemma~\ref{LMefivhid} goes through.
\end{uwgi}

\begin{uwgi}
	Pairs of the form $\pair X f$, where $\map fXX$ is a self-map, are sometimes called discrete dynamical systems, while in universal algebra they are simply \emph{monounary algebras}. It is well-known that for every nonempty set $S$ there exists a free monounary algebra over $S$ and it can be defined as $\pair {S \times \nat} {\un_S}$, where $\un_S \pair s n = \pair s {n+1}$ for every $s \in S$, $\ntr$. Nevertheless, the proof of Lemma~\ref{LMefivhid} is short and straightforward, actually exhibiting the fact that $\pair {S \times \nat} {\un_S}$ is free.	
\end{uwgi}

\section{Main result}

As we have seen, the surjectively universal operator is just a lifting to $\ell_1$ of a universal mapping of a countable infinite set.
The assignment $S \mapsto \ell_1(S)$ is indeed functorial and it is a special case of so-called \emph{free constructions} or, formally, a left adjoint to the forgetful functor assigning the unit ball to a Banach space. We postpone more explanation to the next section.
For the moment, let us formulate a somewhat informal version of our result involving monoid actions.

\begin{tw}\label{THMmejnJeden}
	Let $\M$ be a countable monoid and let $\F \nat$ denote the free object over the set of natural numbers $\nat$. Then there exists an $\M$-action $\uno$ on $\F \nat$ with the following property:
	\begin{enumerate}
		\item[{\rm(U)}] Given a countably generated object $A$, every $\M$-action on $A$ lifts to $\uno$ via a nice epimorphism, namely, there exists a nice epimorphism $\map p {\F \nat} A$ commuting with the $\M$-actions.
	\end{enumerate}
\end{tw}

A little bit of explanation is certainly needed here. By a \emph{free object} we mean an object of some category (e.g. a Banach space, a group, a compact Hausdorff space, etc., see~\cite{BaKa} for more examples) that ``behaves'' like the well-known and classical free group in the sense that every mapping of the generating set $\nat$ into an arbitrary object of the category under consideration extends uniquely to a morphism defined on $\F \nat$. Formal statement involves the associated forgetful functor, which in some cases is not completely obvious. For example, the standard forgetful functor on the category of Banach spaces with non-expansive operators assigns to a Banach space its unit ball instead of the whole space. Being \emph{countably generated} means there exists a surjective mapping from a countable set whose unique extension is a nice epimorphism. As it happens, in some categories (e.g. Banach spaces) epimorphisms are not necessarily surjective, therefore we propose an abstract concept of a \emph{nice} epimorphism that should encode those epimorphisms that we would like to have, for example, quotient mappings between topological spaces or epimorphisms having a right inverse.

Before going into details and presenting the formal statement with a proof, let us state a corollary that extends the Darji -- Matheron result~\cite{DaMa}.

\begin{wn}
	Let $\M$ be a countable monoid. There exists an $\M$-action $\phi$ on $\ell_1$, such that every $\M$-action on a separable Banach space lifts to $\phi$, via a surjective non-expansive linear operator.
\end{wn}

In particular, taking $\M = \Nat$ we obtain a single non-expansive operator, while taking the action of $\Zee_2$, we obtain a surjectively universal isometric involution. Taking the action of $\Zee$, we obtain a surjectively universal linear isometry of $\ell_1$. Finally, taking the free group with countably many generators, we obtain the most complicated surjectively universal action of a countable group on $\ell_1$.

Our method of obtaining surjectively universal actions heavily uses the left adjoint to the forgetful functor, obtaining a concrete description based on the universal actions on sets.
Namely, we prove the following result, where $\fC^\M$ denotes the category of all $\M$-actions on $\fC$, in other words, all covariant functors from $\M$ into $\fC$. A slightly more precise formulation is Theorem~\ref{THMCztyryPiync} below.

\begin{tw}\label{THMkindOfstjupeedMaybe}
	Let $\fC$ be a concrete category with a left adjoint (free) functor $\F$ and let $\M$ be a monoid. Then there exists a functor $\map{\MF}{\fC^\M}{\Ens}$ that is left adjoint to the forgetful functor.
\end{tw}

The forgetful functor from $\fC^\M$ into $\Ens$ simply forgets the $\M$-action and the ``$\fC$-structure''.
Summarizing, once we have a surjectively universal $\M$-action on, say, the set of natural numbers $\nat$, we can lift it using the free functor $\F$, obtaining an $\M$-action on $\F \nat$. This action turns out to be universal in the sense of Theorem~\ref{THMmejnJeden}.

\section{Concrete categories and left adjoints}

We now provide the necessary definitions concerning concrete categories, also introducing the axiomatic notion of \emph{nice} epimorphisms. A good reference concerning concrete categories is~\cite{JoyCats}.
Let us just recall that a category is \emph{concrete} if it admits a faithful (also called: \emph{forgetful}) functor into the category of sets. A functor is \emph{faithful} if it is injective on hom-sets (sets of arrows from a fixed domain to a fixed codomain).

In this section we assume that $\fC$ is a concrete category with a distinguished forgetful functor $\frgc$. That is, every $\fC$-object has its \emph{universe} $\frg{C}$ which is a set, and every $\fC$-morphism $\map f C D$ induces a set mapping $\map{\frg{f}}{\frg{C}}{\frg{D}}$. Furthermore, distinct morphisms between fixed pairs of $\fC$-objects induce distinct mappings between their universes.
In this setting, $\fC$-objects could be regarded as ``structured'' sets and $\fC$-morphisms could be viewed as set mappings preserving those abstract structures.

Denote by $\Ens$ the category of nonempty\footnote{We exclude the empty set just because there are no maps from nonempty sets to $\emptyset$.} sets with arbitrary mappings. Assume $\map{\F}{\Ens}{\fC}$ is a left adjoint to $\frgc$. That is, for each nonempty set $X$, we have a mapping $\map{\eta_X}{X}{\frg{\F X}}$ such that given any $\fC$-object $A$, given a mapping $\map f X \frg A$ there is a unique $\fC$-arrow $\map{\uv f}{\F X}{A}$ satisfying
$$\frg{\uv f} \cmp \eta_X = f.$$
In that case $\F X$ (more precisely, the pair $\pair{\F X}{\eta_X}$) is called the \emph{free $\fC$-object} over $X$.
In fact, if for every nonempty set $X$ one can find a pair $\pair{\F X}{\eta_X}$ as above, then $\F$ extends uniquely to a functor from $\Ens$ to $\fC$ that is left adjoint to $\frgc$. Furthermore, $\eta$ is a natural transformation from the identity of $\Ens$ to $\frg \F$, as visualized in the diagram below.
$$\begin{tikzcd}
	{\frg{\F X}} \ar[r, "{\frg {\F f}}"] & {\frg{\F Y}} \\
	X \ar[u,"\eta_X"] \ar[r, "f"] & Y \ar[u, "\eta_Y"']
\end{tikzcd}$$
Below we note the following well-known properties of the ``$\uv{\;\;\;}$'' operator involved in the construction above.

\begin{lm}\label{LMfudrzibtibfi}
	Assume $X,Y$ are nonempty sets, $A$ is a $\fC$-object, and $\map f X Y$, $\map g Y {\frg A}$ are mappings. Then
	$$\uv{g \cmp f} = \uv g \cmp {\F f}.$$
\end{lm}

\begin{pf} Note that
	$$\frg{\uv g \cmp \F f} \cmp \eta_X = \frg{\uv g} \cmp \frg{\F f} \cmp \eta_X = \frg{\uv g} \cmp \eta_Y \cmp f = g \cmp f.$$
	We have used the fact that $\eta$ is a natural transformation from the identity of $\Ens$ to $\frg{\F(\argum)}$, that is, $\frg{\F f} \cmp \eta_X = \eta_Y \cmp f$.	
\end{pf}

\begin{lm}\label{LMbitgierig}
	Let $\map f X {\frg A}$ be a mapping, where $X \nnempty$ and let $\map g A B$ be a $\fC$-morphism. Then
	$$\uv{\frg g \cmp f} = g \cmp \uv f.$$
\end{lm}

\begin{pf}
	We have	$\frg{g \cmp \uv f} \cmp \eta_X = \frg g \cmp \frg{\uv f} \cmp \eta_X = \frg g \cmp f.$
\end{pf}

\subsection{Nice epimorphisms}

The authors of~\cite{BaKa} define the notion of a \emph{universal free object} (UFO) to be any free object $\F X$ with $X$ an infinite set with the following property: Given an arbitrary $\fC$-object $A$ there is a mapping $\map \iota X A$ such that $\uv \iota$ is an epimorphism.
Obviously, in practice this is meaningful once the category $\fC$ is restricted to objects of a bounded ``size'', ``weight'', ``density'', or another cardinality measure. In that case $X$ can be taken to be of that maximal cardinality. Note also that if $\F X$ is universally free and the cardinality of $X$ is not greater than the cardinality of $Y$ then $\F Y$ is universally free. Indeed, any surjection from $Y$ onto $X$ is right-invertible and hence so is its image under $\F$.

On the other hand, there is a small problem with the notion of an epimorphism. Namely, $f$ is an epimorphism (briefly: epi) if it is left-cancellative, that is
$$g_0 \cmp f = g_1 \cmp f \implies g_0 = g_1,$$
for every compatible morphisms $g_0$, $g_1$.
As it happens, in the category of Banach spaces with non-expansive linear operators (or with all bounded operators) epimorphisms are those mappings whose image is dense. On the other hand, the relevant epimorphisms are surjective operators, due to Banach's open mapping principle.

There are several stronger variants of being an epi, e.g., regular epi, strong epi, extremal epi, etc., see~\cite{JoyCats}. We prefer not to bother with these notions, introducing axiomatically the concept of a \emph{nice} epimorphism.
Namely, we assume that our fixed category $\fC$ has a distinguished 
subcategory $\nepis{\fC}$ satisfying the following axioms.
\begin{enumerate}[itemsep=0pt]
	\item[(N0)] All $\nepis{\fC}$-arrows are epimorphisms.
	\item[(N1)] Right-invertible $\fC$-arrows are in $\nepis{\fC}$.
	\item[(N2)] If $g \cmp h \in \nepis{\fC}$ then $g \in \nepis{\fC}$.
\end{enumerate}
As mentioned above, the $\nepis{\fC}$-arrows will be called \emph{nice}. Note that one can take  $\nepis{\fC}$ to be just the class of all right-invertible $\fC$-arrows. This is the ``minimal'' possibility.
It is easy to check that each of the classes like regular, strong, extremal, and finally all epimorphisms satisfies (N0)--(N2).

\subsection{Generating the objects}

From now on, we assume that $\fC$ is a concrete category with a distinguished class of nice epimorphisms. 
We would like to prevent $\fC$ from any cardinality restrictions, therefore our next task is to replace the concept of a UFO by a more direct notion of being generated by a concrete set.

\begin{df}
	Let $A$ be a $\fC$-object and let $\map f X \frg{A}$ be a map. We say that $f$ \emph{generates} $A$ if $\map {\uv f}{\F X}A$ is a nice epi. If $X \subs \frg A$ and $f$ is the inclusion, then we say that $A$ is \emph{generated} by $X$. We say that $A$ is \emph{$\kappa$-generated} if it is generated by a set of cardinality $<\kappa$. Finally, we say that $A$ is \emph{countably generated} if it is $\aleph_1$-generated, namely, generated by a countable set.
\end{df}

We might also add to the following to the definition above: A $\fC$-object $A$ is \emph{inaccessible} if no mapping into $\frg A$ generates $A$. 

\begin{ex}
	Let $\fC$ be a concrete category such that nice epimorphisms are precisely the right-invertible ones. Then inaccessible objects are those $A \in \ob{\fC}$ such that there is no right-invertible morphism from any $\F X$ to $A$. Such objects do exist, for example, in the case of groups, Banach spaces, and many other concrete categories, as long as we declare nice := right-invertible.
\end{ex}

One would like to reduce generating by a map to generating by a set. This is indeed possible.

\begin{lm}\label{LMstillGenerates}
	Let $\map f X {\frg A}$ be a mapping such that $f$ generates $A$, and let $S$ be such that $\img f X \subs S \subs \frg A$. Then $S$ generates $A$.
\end{lm}

\begin{pf}
	The assumptions lead to the equation $f = e \cmp f_S$, where $f_S$ is the same as $f$, but with co-domain restricted to $S$, and $e$ is the inclusion $S \subs \frg{A}$.
	Lemma~\ref{LMfudrzibtibfi} gives
	$\uv f = \uv e \cmp \F f_S$. Finally, $\uv e$ is nice by (N2), because $\uv f$ was assumed to be nice.
\end{pf}

Below we see that the abstract concept of ``being generated by a set'' behaves in a reasonable way, namely, if $S$ generates $A$ and $S \subs S' \subs \frg A$ then $S'$ generates $A$.

\section{Monoid actions}\label{SectCztyryBilou}

Fix a category $\fC$ and fix a monoid $\M$. We shall consider the category $\fC^\M$ whose objects are pairs of the form $\pair A \phi$, where $A \in \ob{\fC}$ and $\phi$ is an \emph{action} of $\M$ on $A$, that is, a homomorphism from $\M$ into the endomorphism monoid $\End A$. Morphisms are simply those $\fC$-arrows that preserve the $\M$-actions.

Actually, a monoid is just a category with a single object and monoid actions are simply covariant functors from the monoid. Morphisms are natural transformations.

We shall write $\action \M \phi A$ meaning $\map \phi \M {\End A}$ is a monoid action. Given $m \in \M$ we shall write $m^\phi$ for the endomorphism $\phi(m)$. Thus $1^\phi = \id A$ and $(m_0 \cdot m_1)^\phi = m_0^\phi \cmp m_1^\phi$ for every $m_0, m_1 \in \M$.
In the category of sets, we additionally have expressions of the form $m^\phi \cdot x$, namely, the value of $x \in A$ under the mapping $m^\phi = \phi(m)$.

A morphism $\map g A B$ is an \emph{$\M$-homomorphism} from $\action \M \phi A$ to $\action \M \psi B$ if $g \cmp m^\phi = m^\psi \cmp g$ for every $m \in \M$, as in the diagram below.
$$\begin{tikzcd}
	A \ar[r, "m^\phi"] \ar[d,"g"'] & A \ar[d, "g"] \\
	B \ar[r, "m^\psi"'] & B
\end{tikzcd}$$
In other words, $g$ is a natural transformation from $\phi$ to $\psi$.
Note that in case $\fC$ is a concrete category with a forgetful functor into $\Ens$, so is the category of $\M$-actions. Namely, an $\M$-action $\action{\M}{\phi}{A}$ is mapped to an $\M$-action $\action{\M}{\frg{\phi}}{\frg{A}}$ defined in the obvious way: $\frg{\phi}(m) = \frg{\phi(m)}$ for every $m \in \M$. This is in fact the composition of $\phi$ with the forgetful functor.

The free (left adjoint) functor $\F$ from the category of sets to $\fC$ also naturally extends to $\fC^\M$. Indeed, given a set $S$ and a homomorphism $\map \psi \M {\End S}$, we have $\F \psi := \F \cmp \psi$, a homomorphism from $\M$ to $\End({\F S})$. In other words, $\F \psi$ is the $\M$-action on $\F S$ induced by $\phi$.

\begin{lm}\label{LMrbergihweogw}
	Let $\M$ be a monoid, let $S$ be a nonempty set and let $\M \times S$ be endowed with the $\M$-action $\zeta_S$ defined by $m^{\zeta_S} \pair a s = \pair {m \cdot a} s$.
	Then every $\M$-action on $S$ lifts to $\zeta_S$ via a surjection.
\end{lm}

\begin{pf}
	Fix an $\M$-action $\action{\M}{\psi}{S}$.
	Define $\map{q}{\M \times S}{S}$ by $q \pair a s = a^\psi s$.
	We have
	$$q (m^{\zeta_S} \pair a s) = q \pair {m \cdot a} s = (m \cdot a)^\psi s = m^\psi(a^\psi s) = m^\psi(q \pair a s).$$
	Hence $q$ is a homomorphism from $\action{\M}{\zeta_S}{(\M \times S)}$ onto $\action{\M}{\psi}{S}$.
\end{pf}

The $\M$-action described above actually provides a functor from the category $\Ens$ of nonempty sets to the category $\Ens^\M$ of $\M$-actions on nonempty sets. Formally, $\zeta S := \pair{\M\times S}{\zeta_S}$ and, given $\map f ST$, the homomorphism $\map{\zeta f}{\zeta S}{\zeta T}$ is defined by $\zeta f \pair a s := \pair a{fs}$.
We have
$$\zeta f \cmp m^{\zeta_S} = m^{\zeta_T} \cmp \zeta f$$
for every $m \in \M$, therefore $\zeta f$ is indeed a homomorphism of $\M$-actions.
Clearly, $\zeta (f \cmp g) = \zeta f \cmp \zeta g$ and $\zeta \id{S} = \id{\zeta S}$, therefore $\zeta$ is a functor.

The category $\Ens^\M$ of $\M$-actions on nonempty sets has a natural forgetful functor $\frg{\cdot}$, namely, forgetting the action. Obviously, the same applies to any category in place of $\Ens$.

\begin{lm}\label{LMrytilhotssa}
	Let $A \in \ob{\fC}$, let $S$ be a nonempty set. Let $\action{\M}{\psi}{S}$ and $\action{\M}{\phi}{A}$ be $\M$-actions and let $\map f S {\frg A}$ be an $\M$-homomorphism, where the $\M$-action on $\frg A$ is $\frg{\phi}$. Then $\map{\uv{f}}{\F S}{A}$ is an $\M$-homomorphism from $\action{\M}{\F \psi}{\F S}$ to $\action{\M}{\phi}{A}$. 
\end{lm}

\begin{pf}
	Fix $m \in \M$. Then $\map {m^\phi} A A$ and we have
	$$f \cmp m^\psi = \frg{m^\phi} \cmp f,$$
	because $f$ preserves the $\M$-actions.
	Using Lemmas~\ref{LMfudrzibtibfi} and~\ref{LMbitgierig}, we get
	$$\uv f \cmp \F m^\psi = \uv{f \cmp m^\psi} = \uv{\frg{m^\phi} \cmp f} = m^\phi \cmp \uv f.$$
	This completes the proof, because $\F m^\psi = m^{\F \psi}$.
\end{pf}

\begin{tw}\label{THMgdfiuew}
	The functor $\zeta$ defined above is left adjoint to the canonical forgetful functor from $\Ens^\M$ to $\Ens$.
\end{tw}

\begin{pf}
	Given a nonempty set $S$, let $\eta_S$ be the mapping given by $\eta_S(s) = \pair 1 s$.
	Fix a mapping $\map f S \frg X$, where $X$ is a nonempty set endowed with an $\M$-action $\phi$ and $\frg{X}$ is the set $X$ in which the action has been forgotten.
	Define $\uv f$ by setting
	\begin{equation}
		\uv f \pair a s = a^\phi \cdot f(s).
		\label{EQerreg}		
	\end{equation}
	Given $\pair a s \in \M \times S$ and $m \in \M$, we have
	$$m^\phi \cdot \uv f \pair a s = m^\phi \cdot a^\phi \cdot f(s) = (m \cdot a)^\phi \cdot f(s) = \uv f \pair {m a} {f(s)} = \uv f m^\zeta \pair a s,$$
	therefore $\uv f$ is a morphism from $\action{\M}{\zeta}{(\M \times S)}$ to $\action{\M}{\phi}{S}$. Its uniqueness is obvious, as formula (\ref{EQerreg}) is the only way of defining an $\M$-action coherent with $f$.
\end{pf}


Recall that an object $A$ is \emph{$\kappa$-generated} if it is generated by a set of cardinality $\loe \kappa$.

\begin{tw}\label{THMCztyryCztyryHa}
	Let $\M$ be a monoid of cardinality $\loe \kappa$, where $\kappa$ is infinite. Let $\fC$ be a concrete category with nice epimorphisms. Assume $\F$ is a left adjoint to the forgetful functor.
	Then $\F \kappa$ has an $\M$-action $\uno = \F \zeta_\kappa$ that lifts, via nice epimorphisms, all $\M$-actions on $\kappa$-generated $\fC$-objects.
\end{tw}

\begin{pf}
	We first show how the $\M$-action on a $\fC$-object lifts nicely to an $\M$-action induced by sets.
	Fix a $\fC$-object $A$ generated by a set $S \subs \frg A$ of cardinality $\loe \kappa$ and fix an $\M$-action $\action{\M}{\phi}{A}$.
	We may assume that $S$ is $\M$-invariant, simply by closing it off with respect to all the $\M$-orbits. Since the cardinality of $\M$ is $\loe \kappa$, the cardinality of $S$ remains $\loe \kappa$ and by Lemma~\ref{LMstillGenerates}, $S$ still generates $A$.
	Now, we have an $\M$-action
	$\action{\M}{\psi}{S}$, namely, $\psi$ is the restriction of $\frg \phi$ to $S$.
	In other words, $m^\psi = \frg{m^\phi} \rest S$ for every $m \in \M$.
	Let $\iota$ be the inclusion $S \subs \frg A$. Note that $\iota$ is an $\M$-homomorphism, therefore by Lemma~\ref{LMrytilhotssa} $\uv \iota$ is an $\M$-homomorphism from the $\M$-action $\F \psi = \F \cmp \psi$ on $\F S$ onto the $\M$-action $\phi$.
	Furthermore, $\uv \iota$ is a nice epi, because $S$ generates $A$. Thus, $\phi$ nicely lifts to $\F \psi$.
	Finally, Lemma~\ref{LMrbergihweogw} shows that $\F \psi$ lifts to $\F \zeta_\kappa$ via a nice epimorphism, because every surjection of sets is right-invertible.
\end{pf}

We are left with the proof of Theorem~\ref{THMkindOfstjupeedMaybe}.
We have a functor $\map{\F^\M}{\Ens^\M}{\fC^\M}$ defined by
$$(\F^\M) \pair X \phi := \pair{\F X}{\F \phi} \oraz (\F^\M) f := \F f$$
for $f \in {\Ens^\M}$. Finally:

\begin{tw}\label{THMCztyryPiync}
	The functor
	$$\MF := \F^\M \cmp \zeta.$$
	is left adjoint to the canonical forgetful functor from $\fC^\M$ to $\Ens$.
\end{tw}

\begin{pf}
	Fix $\pair A \phi$ in $\fC^\M$. Note that $\frg{\pair A \phi} = \frg A$. Fix a mapping $\map{f}{S}{\frg A}$, where $S$ is a nonempty set.
	Define $\map{\tilde{f}}{\M \times S}{\frg A}$ by $\tilde{f}\pair x s := x^{\frg{\phi}} \cdot f(s)$, where $\frg{\phi}$ is the composition of $\phi$ with the forgetful functor. Note that this is the unique possibility of defining $\tilde f$ extending $f$ and preserving the $\M$-actions.
	Let $\uv f$ be the unique extension of $\tilde f$ to $\F (\M\times S)$. By Lemma~\ref{LMrytilhotssa}, $\uv f$ is an $\M$-homomorphism.
\end{pf}

\paragraph{Conclusion.}
The story presented above exhibits the following phenomenon. Namely, given a concrete category $\fC$ and a monoid $\M$, we have the following commutative diagram
$$\begin{tikzcd}
	\Ens \ar[rr, "\F"] \ar[dd, "\zeta"] \ar[ddrr, "\MF", bend right] & & \fC \ar[ll, bend right, dotted]\\
	& & \\
	{\Ens^\M} \ar[rr, "\F^\M"'] \ar[uu, bend left, dotted] & & {\fC^\M} \ar[uu, bend right, dotted] \ar[ll, bend left, dotted]
\end{tikzcd}$$
where dotted arrows indicate the forgetful functors.
This applies to every concrete category, in particular, to all categories mentioned in~\cite{BaKa}. Note that the concept of nice epimorphisms is needed only for obtaining a better universality result in Theorem~\ref{THMCztyryCztyryHa} above. On the other hand, Theorem~\ref{THMCztyryPiync} is valid for every concrete category admitting a left adjoint to its forgetful functor.

Finally, the category $\Ens$ of nonempty sets can be replaced by another category $\fK$ in which there are free monoid actions taking the role of $\zeta$. In that case, $\F$ is a left adjoint to a fixed faithful functor from $\fC$ to $\fK$.

\section{Remarks on surjectively universal mappings}

A surjectively universal mapping of a countable infinite set is of course not unique, although the one we have proposed just before Lemma~\ref{LMefivhid} is the canonical one, namely, it can be viewed as the free monounary algebra over a countable infinite set. The mapping proposed by Darji and Matheron~\cite{DaMa} (repeated in~\cite{BaKa}) is actually injectively universal. It lifts maps that have at least one fixed point, that is why the free operator on $\ell_1$ induced by this mapping is actually surjectively universal. Our goal here is to explain the reasons why the proof of~\cite{DaMa} works, while the much more general argument of Balcerzak and Kania~\cite{BaKa} works only in some cases, namely, when it is enough to restrict attention to mappings that have at least one fixed point.

Given a self-mapping $\map f XX$ its associated \emph{digraph} has $X$ as the set of vertices and the directed edges are of the form $\pair{x}{fx}$, where $x \in X$. Note that every digraph whose out-degree is precisely one encodes a self-mapping of its vertex set. Of course, loops are allowed---they simply encode the fixed points.
Let us denote by $\Gr f$ the graph described above. Note that in fact $\Gr f$ is just the graph of $f$, viewed as a binary relation. In any case, we can now talk about the components of $f$, having in mind the associated digraph.

We say that a digraph $G$ is \emph{natural} if there is a homomorphism from $G$ onto $\Nat$, where each $n \in \Nat$ is connected to its immediate successor $n+1$. In particular, a natural digraph has no loops. The natural digraph of $\Nat$ corresponds to the shift $n \mapsto n+1$.

In what follows, we restrict attention to self-mappings of countable sets, the general case can be easily deduced from our considerations.

\begin{tw}\label{THMxawasy}
	Assume $\map w X X$ is a self-mapping of a countable infinite set. Then $w$ is surjectively universal if and only if it satisfies the following conditions.	\begin{enumerate}[itemsep=0pt]
		\item[{\rm(I)}] $\Gr w$ has infinitely many components.
		\item[{\rm(W)}] Each component of $\Gr w$ is natural.
	\end{enumerate}
\end{tw}

We could have reformulated (W) by saying that the whole graph $\Gr w$ is natural, however we would like to stress the fact that each component is infinite, which of course follows from the definition of being natural.

\begin{pf}
	Clearly, $w$ is surjectively universal if and only if it lifts $\un$.
	
	Suppose $w$ lifts $\un$, namely, $p \cmp w = \un \cmp p$, where $p$ is a surjection. Note that $p$ preserves the digraph structures therefore $\Gr w$ must have infinitely many components, because $\un$ has so.
	Now recall that $\map{\un}{\omega^2}{\omega^2}$ is defined by $\un \pair mn = \pair {m+1}n$. Given $x \in X$ let $r(x) = m$ if $p(x) = \pair mn$. Observe that $r(y) = r(x) + 1$ whenever $y = w x$. Thus, $r$ witnesses that $\Gr w$ is natural---its image might not be onto, but we may shift it to the left appropriately.
	
	Conversely, suppose $\Gr f$ satisfies (I) and (W).
	Let $\sett{C_n}{\ntr}$ enumerate all components of $\Gr f$.
	For each component $C_n$ fix a witness $\map {p_n}{C_n}\Nat$ for being natural.
	Given $x \in C_n$, let $p(x) = \pair {p_n(x)}n$.
	Since each $p_n$ is onto, $p$ is a surjection.
	Finally, if $y = fx$ and $x \in C_n$ then $p_n(y) = p_n(x) + 1$, therefore $\un \cmp p = p \cmp f$, showing that $f$ lifts $\nu$.
\end{pf}

In particular, a surjectively universal map cannot have $\Zee$-orbits, i.e., orbits of the form $\setof{f^n x}{n \in \Zee}$.
On the other hand, we have the following fact, explaining why the operator constructed in~\cite{DaMa} is projectively universal.

\begin{prop}
	Assume $\map w XX$ is a self-map of a countable set with infinitely many natural components.
	Then $w$ lifts every countable self-map with a fixed point.
\end{prop}

\begin{pf}
	Let $N \subs X$ be the union of all natural components of $w$. Fix $\map f S S$ with $S$ countable such that $f(s_0) = s_0$. By Theorem~\ref{THMxawasy}, $w \rest N$ lifts $f \rest S \setminus \sn{s_0}$. Mapping $X \setminus N$ to $s_0$, we obtain a lifting of $f$ to $w$. In case $N = X$, we leave out one of the infinitely many natural components of $w$ and map it onto $s_0$.
\end{pf}

\paragraph{Acknowledgments.} The author would like to thank the participants of the Set Theory Seminar at Pavol Jozef Šafárik University in Košice, led by Jaroslav Šupina, where a preliminary version of this note was presented. Special thanks are due to Miroslav Plo\v s\v cica, for pointing out the use of free mono-unary algebras. The author is also grateful to Adam Barto\v s, Tristan Bice, and Tomasz Kania for useful remarks.

\end{document}